\documentclass[11pt]{amsart}
\usepackage{amsmath,amsthm,amsfonts,amscd,amssymb,eucal,latexsym}

\theoremstyle{definition}
\newtheorem{theorem}{Theorem}[section]
\newtheorem{corollary}[theorem]{Corollary}
\newtheorem{lemma}[theorem]{Lemma}

\newtheorem{remark}[theorem]{Remark}

\newtheorem{definition}[theorem]{Definition}

\newcommand{\topol}{_{\text{\rm top}}}

\newcommand{\CPA}{{\rm CPA}}
\newcommand{\rcp}{{\rm rcp}}

\newcommand{\sa}{{\rm sa}}
\newcommand{\Aff}{{\rm Aff}}
\newcommand{\sep}{{\rm sep}}

\newcommand{\spn}{{\rm span}}

\newcommand{\interior}{{\rm int}}

\allowdisplaybreaks

\begin{document}

\title[Positive topological entropy and $\ell_1$]{Positive topological
entropy and $\ell_1$}

\author{David Kerr}
\author{Hanfeng Li}
\address{\hskip-\parindent
David Kerr \\Dipartimento di Matematica, Universit\`{a} di Roma
``La Sapienza,'' P.le Aldo Moro, 2, 00185 Rome, Italy}
\email{kerr@mat.uniroma1.it}
\address{\hskip-\parindent
Hanfeng Li \\Department of Mathematics, University of Toronto, Toronto,
Ontario M5S 3G3, Canada}
\email{hli@fields.toronto.edu}
\date{March 13, 2003}

\begin{abstract}
We characterize positive topological entropy for quasi-state space
homeomorphisms induced from $C^*$-algebra automorphisms in terms of
dynamically generated subspaces isomorphic to $\ell_1$. This geometric
condition is also used to give a description of the topological Pinsker
algebra. In particular we obtain a geometric characterization of positive 
entropy for topological dynamical systems,
as well as an analogue for completely positive topological entropy
of Glasner and Weiss's combinatorial characterization of completely positive
Kolmogorov-Sinai entropy.
\end{abstract}

\maketitle

\section{Introduction}

In \cite{GW} E. Glasner and B. Weiss showed that if
a homeomorphism from a compact metric space $X$ to itself has zero
topological entropy, then so does the induced homeomorphism on the space
of probability measures on $X$ with the weak$^*$ topology. One of the
two proofs they gave of this striking result established a remarkable
connection between topological dynamics and the local
theory of Banach spaces. The key geometric fact is the exponential
dependence of $k$ on $n$ given an
approximately isometric embedding of $\ell^n_1$ into $\ell^k_\infty$,
which they deduced from the work of T. Figiel, J. Lindenstrauss, and
V. D. Milman on almost Hilbertian sections of unit balls in Banach
spaces \cite{FLM}.

The first author showed in \cite{EID} that Glasner and Weiss's
geometric approach can be conceptually simplified from a
functional-analytic viewpoint using Voiculescu-Brown entropy
(see Remark 3.10 in \cite{EID}) and also
more generally applied to show that if an automorphism of a separable
exact $C^*$-algebra has zero Voiculescu-Brown entropy then the induced
homeomorphism on the quasi-state space has zero topological entropy.
In this case the crucial Banach space fact is the exponential dependence of
$k$ on $n$ given an approximately isometric embedding of $\ell^n_1$ into
the matrix $C^*$-algebra $M_k$ \cite[Lemma 3.1]{EID}, which can be deduced
from the work of N. Tomczak-Jaegermann on the Rademacher type $2$ constants
of Schatten $p$-classes \cite{TC}.
Following a procedure similar to that of \cite{GW} but in a
matrix setting, the key Banach space map in \cite{EID} is constructed by
passing suitable elements in the $C^*$-algebra through a matrix algebra
$M_k$ in a completely positive approximation and then evaluating on
states which register the entropic growth, thus producing a
map from the trace class
$C^k_1$ to $\ell^n_\infty$ which, after composing with a suitable
projection and taking the dual, yields the desired approximately
isometric embedding of $\ell^n_1$ into $M_k$.

As observed in Remark 3.11 of \cite{EID}, Lemma 3.1 of \cite{EID}
can also be used directly at the $C^*$-algebra level with respect to
the completely positive approximations that appear in the definition
of Voiculescu-Brown entropy. Indeed if a finite set $\Omega$ of elements
in the $C^*$-algebra forms a standard basis for an isomorphic copy
of $\ell^n_1$ then a completely positive approximation to $\Omega$ of
suitably fine tolerance will yield an isomorphic embedding of $\ell^n_1$
into the associated matrix algebra with a similar isomorphism constant
(see Lemma~\ref{L-ell1} in Section~\ref{S-tPa}). We may thus ask
if the geometric phenomena in \cite{GW,EID} connected to positive 
topological entropy on the quasi-state space can be entirely captured 
within the Banach space structure of the $C^*$-algebra itself.
The principal aim of this note is to answer this question in the affirmative.
Specifically, we show that positive entropy on the quasi-state space
can be characterized by the existence of a subspace of the $C^*$-algebra
isomorphic to $\ell_1$ which is generated in a canonical way by a single
element along a set of iterates of positive density. This forms the
content of Section~\ref{S-qs}, the first of the two sections comprising the 
main body of the paper. 

In Section~\ref{S-tPa}
we apply arguments similar to those from Section~\ref{S-qs} in conjunction
with results from \cite{EID} to obtain a geometric characterization of the
topological Pinsker algebra. This dynamical object is the $C^*$-algebraic
manifestation of the maximal zero entropy factor of a topological
dynamical system \cite{ZEF} and is an analogue of the Pinsker
$\sigma$-algebra in ergodic theory. We show that the self-adjoint elements of
the topological Pinsker algebra are precisely those which do not dynamically
generate in a canonical way a subspace isomorphic to $\ell_1$ along
a set of iterates of positive density. As a corollary we obtain
a geometric characterization of positive entropy for topological dynamical
systems which answers a question posed by V. Pestov after the first author's
talk at the Winter 2002 CMS Meeting in Ottawa, as well as
an analogue for completely positive topological entropy
of Glasner and Weiss's combinatorial characterization of completely
positive Kolmogorov-Sinai entropy \cite{GW}.
\medskip

\noindent{\it Acknowledgements.} D. Kerr was supported by the
Natural Sciences and Engineering Research Council of Canada. This
work was carried out during his stay at the University of Rome
``La Sapienza'' over the 2002--2003 academic year. He thanks
Claudia Pinzari at the University of Rome ``La Sapienza'' for her
generous hospitality. H. Li was supported jointly by the
Mathematics Department of the University of Toronto and NSERC
Grant 8864-02 of George A. Elliott.

\section{A geometric characterization of positive entropy on the
quasi-state space}\label{S-qs}

We begin by establishing notation. Let $A$ be a $C^*$-algebra and
$\alpha$ an automorphism of $A$. Following the notation of
\cite{EID}, we denote by $S_\alpha$ the homeomorphism of the
closed unit ball of $A^*$ (with the weak$^*$ topology) given by
$S_{\alpha}(\sigma ) = \sigma\circ\alpha$, and by
$\tilde{T}_\alpha$ the restriction of $S_\alpha$ to the quasi-state space
$Q(A)$, i.e., the weak$^*$ compact convex
set of positive linear functionals on $A$ of norm at most one. Also,
when $A$ is unital we denote by $T_\alpha$ the restriction of
$\tilde{T}_\alpha$ to the state space $S(A)$, i.e., the weak$^*$
compact convex set of positive unital linear functionals on $A$.
The real linear subspaces of self-adjoint elements of $A$ and
$A^*$ will be denoted by $A_{\sa}$ and $(A^*)_{\sa}$,
respectively.

Let $X$ be a compact topological space and $T : X\to X$ a
homeomorphism. Let $d$ be a pseudo-metric on $X$. A set $E\subseteq X$ is
said to be {\it $(n,\varepsilon )$-separated (with respect to $T$ and $d$)}
if for every $x,y\in E$ with $x\neq y$ there exists a
$0\leq k \leq n-1$ such that $d(T^k x , T^k y) > \varepsilon$.
We denote by $\sep_n (T,\varepsilon )$
the largest cardinality of an $(n,\varepsilon )$-separated set, and we
define
$$ h_d (T) = \sup_{\varepsilon >0} \limsup_{n\to\infty}\frac1n
\log\sep_n (T,\varepsilon ) . $$
If $d$ is a metric compatible with the topology on $X$ then
$h_d (T)$ agrees with the topological entropy $h_{\topol}(T)$ defined
as the supremum of $h_{\topol} (T,\mathcal{U})$
over all finite open covers $\mathcal{U}$ of $X$, where
$$ h_{\topol} (T,\mathcal{U}) = \lim_{n\to\infty} \frac1n \log
N(\mathcal{U}\vee T^{-1}\mathcal{U}\vee\cdots\vee T^{-(n-1)} \mathcal{U}) $$
and $N(\cdot )$ denotes the smallest cardinality of a subcover.
See \cite{DGS,Wal} for references on topological entropy.

We will have occasion to consider Banach spaces over both the real and
complex numbers. If the scalar field is not clear from the context
then it will be assumed to be the complex numbers unless the notation
is tagged with an $\mathbb{R}$. In particular the Banach spaces $\ell_1$
over some index set which appear in the statements of our results are
complex, although this is not essential as the analogous
statements using real scalars are also valid.

\begin{definition}
Let $X$ be a complex Banach space and $\alpha : X\to X$ an
isometric isomorphism. Let $x\in X$. We say that an infinite subset
$I\subseteq\mathbb{Z}$ is an {\it $\ell_1$ isomorphism set for $x$} if
there is an isomorphism from $\ell^I_1$ to
$\overline{\spn} \{ \alpha^i (x) : i\in I \}$ sending
the standard basis element of $\ell^I_1$ associated with
$i\in I$ to $\alpha^i (x)$.
\end{definition}

Recall that the {\it density} of a set $J\subseteq\mathbb{Z}$ is defined as
the limit
$$ \lim_{n\to\infty} \frac{| J \cap \{ -n, -n+1 , \dots , n \} |}{2n+1} $$
if it exists. If $X$ and $Y$ are Banach spaces and $\Gamma : X \to
Y$ is an isomorphism then we say that $\Gamma$ is a {\it
$K$-isomorphism} if $\| \Gamma \| \| \Gamma^{-1} \| \leq K$.

\begin{theorem}\label{T-qs}
Let $A$ be a separable $C^*$-algebra and $\alpha$ an automorphism
of $A$. Then the following are equivalent:
\begin{enumerate}
\item $h_{\topol} (\tilde{T}_\alpha ) = \infty$,
\item there exist an $a\in A$, constants $K, d>0$, a sequence
$\{n_k\}_{k\in \mathbb{N}}$ in $\mathbb{N}$ going to $\infty$, and sets
$I_k\subseteq \{0, 1, \cdots, n_k-1\}$ of cardinality at least
$dn_k$ for each $k\in\mathbb{N}$ such that, for each $k\in\mathbb{N}$,
the linear map from $\ell^{I_k}_1$ to
$\spn \{ \alpha^i (a) : i\in I_k\}$
which sends the standard basis element of $\ell^{I_k}_1$
associated with $i\in I_k$ to $\alpha^i (a)$ is a $K$-isomorphism,
\item there exists an $a\in A$ with an $\ell_1$ isomorphism set of positive
density.
\end{enumerate}
We may moreover take $a$ in (2) and (3) to be in any given
countable total subset $W$ of $A_{\sa}$. When $A$ is unital then
the conditions (1)--(3) are also equivalent to
\begin{enumerate}
\item[(4)] $h_{\topol} (T_\alpha ) = \infty$.
\end{enumerate}
\end{theorem}

\begin{proof}
(1)$\Rightarrow$(2). Let $W = \{ a_1 , a_2 , a_3 , \dots \}$ be
a countable total subset of $A_{\sa}$. Without loss of generality we
may assume that $W$ is a subset of the unit ball $B_1 (A_{\sa} )$. Define
a metric on $B_1 (A^* )$ by
$$ d(\sigma ,\omega ) = \sum_{j=1}^{\infty} 2^{-j}
|(\sigma -\omega )(a_j )| . $$
Then there exist $\varepsilon > 0$, $\lambda > 0$,
a sequence $\{ n_k \}_{k\in\mathbb{N}}$
in $\mathbb{N}$ going to $\infty$, and an $(n_k , 4\varepsilon )$-separated
subset $E_k$ of $Q(A)$ with cardinality at least $e^{\lambda n_k}$ for each
$k\in\mathbb{N}$. Take an $N > 0$ with $2^{-N} < \varepsilon$ and define
the pseudo-metric
$$ d' (\sigma, \omega ) = \sum_{j=1}^N 2^{-j} |(\sigma -\omega )(a_j )| $$
on $B_1 (A^* )$. Then
$| d(\sigma, \omega ) - d' (\sigma, \omega )| < \varepsilon$ for all
$\sigma, \omega\in Q(A)$, and so $E_k$ is
$(n_k, 2\varepsilon )$-separated with respect to $d'$. In
particular for distinct $f,g \in E_k$ there exist $1\leq j\leq N$ and
$0\leq i\leq n_k - 1$ such that $|(f-g)(\alpha^i (a_j ))| > \varepsilon$.

Define next an $\mathbb{R}$-linear map $\phi : (A^*)_{\sa} \to
(\ell^{N\times n_k}_{\infty})_{\mathbb{R}}$ by $(\phi (f))_{ji} =
f(\alpha^i (a_j ))$, where the standard basis of $(\ell^{N\times
n_k}_{\infty})_{\mathbb{R}}$ is indexed by $\{ 1, \dots, N \}
\times \{ 0, \dots, n_k - 1 \}$. Then $\phi (E_k )$ is
$\varepsilon$-separated. As in the proof of Proposition 2.1 of
\cite{GW}, there are constants $d, \delta > 0$ depending only on
$\varepsilon$ and $\lambda$ such that for all sufficiently large
$k$ there is a set $J_k \subseteq \{ 1, \dots, N \} \times \{ 0,
\dots, n_k - 1 \}$ with
\begin{enumerate}
\item[(i)] $|J_k|\geq dNn_k$, and \item[(ii)] $\pi (\phi (B_1
((A^* )_{\sa})))\supseteq B_{\delta}
((\ell^{J_k}_{\infty})_{\mathbb{R}})$, where $\pi : (\ell^{N\times
n_k}_{\infty})_{\mathbb{R}}\to (\ell^{J_k}_{\infty})_{\mathbb{R}}$
is the canonical projection.
\end{enumerate}
Then for any such sufficiently large $k$ there exist some $1\leq
j(k)\leq N$ and a set $I_k\subseteq \{ 0, \dots, n_k - 1 \}$ such
that $|I_k |\geq dn_k$ and $\{ j(k) \} \times I_k \subseteq J_k$.
Consequently $\pi' (\phi (B_1 ((A^* )_{\sa})))\supseteq B_{\delta}
((\ell^{I_k}_{\infty})_{\mathbb{R}})$, where $\pi':(\ell^{N\times
n_k}_{\infty})_{\mathbb{R}}\to (\ell^{I_k}_{\infty})_{\mathbb{R}}$
is the canonical projection. The dual $(\pi' \circ \phi )^*$ is an
injection of $((\ell^{I_k}_{\infty})_{\mathbb{R}})^* =
(\ell^{I_k}_1 )_{\mathbb{R}}$ into $((A^* )_{\sa})^*$ and the norm
of the inverse of this injection is bounded above by
$\delta^{-1}$.  Notice that $A_{\sa} \subseteq ((A^*)_{\sa})^*$,
and from our definition of $\phi$ it is clear that $(\pi' \circ
\phi )^*$ sends the standard basis element of $\ell^{I_k}_1$
associated with $i\in I_k$ to $\alpha^i (a_{j(k)})$.

Let $\Gamma_k$ be the complexification of the map $(\pi' \circ
\phi )^* : (\ell^{I_k}_1 )_{\mathbb{R}}\to
\spn_{\mathbb{R}}\{\alpha^j (a_{j(k)}) : j\in I_k \}\subseteq A$.
Since $a_{j(k)}$ is self-adjoint, the inverse of $\Gamma_k$
evidently has norm no bigger than $K := 2\delta^{-1}$. Since
$1\leq j(k)\leq N$ there is a $j_0$ such that $j(k) = j_0$ for
infinitely many $k$. By taking a subsequence of $\{ n_k \}_{k\in
\mathbb{N}}$ we may assume that $j(k) = j_0$ for all $k$. Now set
$a = a_{j(0)}$.

(2)$\Rightarrow$(1). Multiplying $a$ by a scalar we may assume
that $\| a \| = 1$. Take a dense sequence $a_1 = a, a_2, a_3 , \dots$
in the unit ball $B_1 (A)$ and define a metric on $B_1 (A^* )$ in the same
way as in the first part of the proof.
Denote $\spn \{\alpha^i (a) : i\in I_k \}$ by $V_k$, and let
$\Gamma_k$ denote the linear map from $\ell^{I_k}_1$ to $V_k$ sending
the standard basis element of $\ell^{I_k}_1$ associated with $i\in I_k$
to $\alpha^i (a)$. For each $f\in
(\ell^{I_k}_1 )^*$ we have $(\Gamma^{-1}_k )^* (f)\in (V_k )^*$ and
$\| (\Gamma^{-1}_k )^* (f) \| \leq K \| f \|$. By the Hahn-Banach
theorem we may extend $(\Gamma^{-1}_k )^* (f)$ to an element in
$A^*$, which we will still denote by $(\Gamma^{-1}_k )^* (f)$. Let
$0 < \varepsilon < (2K)^{-1}$, and let $M = \lfloor
(2K\varepsilon)^{-1} \rfloor$ be the largest integer no greater than
$(2K\varepsilon )^{-1}$. Let $\{ g_i : i\in I_k \}$ be the standard basis
of $(\ell^{I_k}_1 )^* = \ell^{I_k}_{\infty}$. For each $f\in \{ 1,
\dots, M \}^{I_k}$ set $\tilde{f} = \sum_{i\in I_k} 2f(i)\varepsilon
g_i$. Then $f' := (\Gamma^{-1}_k )^* (\tilde{f} )$ is in $B_1 (A^* )$.

We claim that the set $\{ f' : f\in \{1, \cdots, M \}^{I_k} \}$ is $(n_k,
\varepsilon )$-separated. Suppose $f, g \in \{1, \cdots, M\}^{I_k}$ and
$f(i) < g(i)$ for some $i\in I_k$. Then
\begin{align*}
d(S^i_{\alpha} (f'), S^i_{\alpha} (g'))
&\geq | (S^i_{\alpha} (f'))(a) - (S^i_{\alpha} (g'))(a) | \\
&= | f'(\alpha^i (a)) - g'(\alpha^i (a)) | \\
&= 2(g(i) - f(i))\varepsilon \\
&> \varepsilon.
\end{align*}
This establishes our claim. Thus $\sep_{n_k} (S_{\alpha} ,
\varepsilon )\geq M^{|I_k|} \geq M^{dn_k}$. It follows that
$h_{\topol} (S_{\alpha} )\geq d\log M$. Letting $\varepsilon\to 0$
we get $h_{\topol} (S_{\alpha} ) = \infty$. Thus $h_{\topol}
(\tilde{T}_{\alpha} ) = \infty$ by Lemma 3.3 of \cite{EID}.

(2)$\Rightarrow$(3). Let $Y_a$ be the collection of sets
$I\subseteq\mathbb{Z}$ such that the linear map from $\ell^I_1$ to
$\overline{\spn} \{ \alpha^i (a) : i\in I \}$ which sends the
standard basis element of $\ell^I_1$ associated with $i\in I$ to
$\alpha^i (a)$ is a $K$-isomorphism. If we identify the subsets of 
$\mathbb{Z}$ with elements of $\{ 0,1 \}^{\mathbb{Z}}$
via their characteristic functions then $Y_a$ is a closed
shift-invariant subset of $\{ 0,1 \}^{\mathbb{Z}}$. It follows by
the argument in the second paragraph of the proof of Theorem 3.2
in \cite{GW} that $Y_a$ has an element $J$ with density at least
$d$. Clearly $J$ is an $\ell_1$ isomorphism set for $a$.

(3)$\Rightarrow$(2). This is immediate in view of the fact that $\alpha$ is
isometric.

Finally, we note that the equivalence (1)$\Leftrightarrow$(4) in
the unital case follows from Lemmas 2.2 and 3.3 of \cite{EID}.
\end{proof}

\begin{remark}\label{R-posentr}
(i) It is easily seen that for the countable subset $W\subseteq A_{\sa}$
we need in fact only assume that $\bigcup_{j\in\mathbb{Z}} \alpha^j (W)$
is total in $A$.

(ii) It is desirable to be able to relax the condition $W\subseteq A_{\sa}$
to $W\subseteq A$. However, we have been unable to determine if this is
possible.

(iii) The proof of Theorem 3.4 in \cite{EID} uses the fact that the
linear maps in the definition of Voiculescu-Brown entropy are
self-adjoint. Using Theorem~\ref{T-qs} above and Lemma 3.1 of
\cite{EID} (or more precisely Lemma~\ref{L-ell1} in
Section~\ref{S-tPa} below) one can prove Theorem 3.4 of \cite{EID}
using only the fact that these linear maps are contractions, and
hence more in the spirit of a Banach space approach.

(iv) The implication (1)$\Rightarrow$(2) is useful for proving
$h_{\topol}(\tilde{T}_{\alpha})=0$ for certain $\alpha$. See
Corollary~\ref{C-free} below for an example.
\end{remark}

\begin{corollary}\label{C-free}
Let $\sigma$ be a permutation of $\mathbb{Z}$. Let $\sigma_*$ be the
corresponding free permutation automorphism of the full free group
$C^*$-algebra $C^* (F_{\mathbb{Z}})$ sending $u_j$ to $u_{\sigma(j)}$, where
$\{u_j\}_{j\in\mathbb{Z}}$ is the set of canonical unitaries of
$C^* (F_{\mathbb{Z}})$. Then $h_{\topol}(T_{\sigma_*}) = \infty$ if and only
if $\sigma$ has an infinite orbit.
\end{corollary}

\begin{proof} The ``if'' part follows from the proof of Proposition
2.4 of \cite{EID}. For the ``only if'' part, suppose that $\sigma$ has no
infinite orbits. For each $g\in F_{\mathbb{Z}}$ let $u_g$ be the
corresponding unitary in $C^*(F_{\mathbb{Z}})$. Then
by taking $W = \{ (u_g + u^*_g )/2 , (u_g - u^*_g )/(2i) : g\in
F_{\mathbb{Z}} \}$ in Theorem~\ref{T-qs} we get
$h_{\topol}(T_{\sigma_*}) = 0$, as desired.
\end{proof}

\section{A geometric characterization of the topological Pinsker
algebra}\label{S-tPa}

In \cite{ZEF} F. Blanchard and Y. Lacroix introduced the maximal
zero entropy factor of a topological dynamical system. This was
called the topological Pinsker factor in \cite{Gl} and is defined
as follows. Let $X$ be a compact metric space and $T : X\to X$ a
homeomorphism. A pair $(x,y)\in X\times X$ with $x\neq y$ is
called an {\it entropy pair} if $h_{\topol} (T,\mathcal{U}) > 0$
for every two-element open cover $\mathcal{U} = \{ U,V \}$ with
$x\in\interior (X\setminus U)$ and $y\in\interior (X\setminus V)$
\cite{Disj}. The {\it topological Pinsker factor} is the quotient
system arising from the closed $T$-invariant equivalence relation
on $X$ generated by the collection of entropy pairs.

Let $\alpha_T$ be the automorphism of $C(X)$ given by
$\alpha_T (f) = f\circ T$ for all $f\in C(X)$. The topological Pinsker factor
corresponds at the $C^*$-algebra level to the $\alpha_T$-invariant
$C^*$-subalgebra $P_{X,T}$ of $C(X)$ whose elements are those functions
$f\in C(X)$ which satisfy $f(x) = f(y)$ for
every entropy pair $(x,y)$. We refer to $P_{X,T}$ as the
{\it topological Pinsker algebra}. It is an analogue of the Pinsker
$\sigma$-algebra in ergodic theory. We also refer to the real Banach algebra
$(P_{X,T} )_{\sa}$ of self-adjoint elements of $P_{X,T}$ as the
{\it real topological Pinsker algebra}. The goal of this section is to obtain
a geometric description of the elements of $(P_{X,T} )_{\sa}$ and
$P_{X,T}$.

We will employ the notation used in Section 2 of \cite{EID} for the
terms involved in the definition of Voiculescu-Brown entropy \cite{Br},
in particular the completely positive rank $\rcp (\Omega , \delta )$ and
the local Voiculescu-Brown entropy $ht(\alpha , \Omega )$. For other
notation and terminology see the beginning of Section~\ref{S-qs}.

The following lemma essentially appeared in Remark 3.11 of
\cite{EID}.

\begin{lemma}\label{L-ell1}
Let $A$ be an exact $C^*$-algebra. Let $a_1 , \dots , a_n \in A$ and suppose
that the linear map $\Gamma : \ell^n_1 \to \spn \{ a_1 , \dots , a_n \}$
sending the $i$th standard basis element of $\ell^n_1$ to $a_i$
for each $i = 1, \dots , n$ is an isomorphism. Let $\delta > 0$ be such
that $\delta < \| \Gamma^{-1} \|^{-1}$. Then
$$ \log\rcp (\{ a_1 , \dots , a_n \} , \delta ) \geq nc \| \Gamma \|^{-2}
(\| \Gamma^{-1} \|^{-1} - \delta )^{2} $$
where $c>0$ is a universal constant.
\end{lemma}

\begin{proof}
Let $\pi : A\to\mathcal{B} (\mathcal{H})$ be a faithful
$^*$-representation, and suppose $(\phi , \psi , B)
\in\CPA (\pi , \{ a_1 , \dots , a_n \} , \delta )$. For any linear
combination $\sum c_i a_i$ of the $a_i$'s we have
\begin{align*}
\Big\| \sum c_i a_i \Big\| &\leq \Big\| \pi\Big( \sum c_i a_i \Big)
- (\psi\circ\phi )\Big( \sum c_i a_i \Big) \Big\|
+ \Big\| (\psi\circ\phi )\Big( \sum c_i a_i \Big) \Big\| \\
&\leq \delta \sum |c_i | + \Big\| \phi \Big( \sum c_i a_i \Big) \Big\| \\
&\leq \delta \| \Gamma^{-1} \|
\Big\| \sum c_i a_i \Big\| + \Big\| \phi \Big( \sum c_i a_i \Big) \Big\|
\end{align*}
and so $\big\| \phi \big( \sum c_i a_i \big) \big\| \geq (1 - \delta
\| \Gamma^{-1} \| ) \big\| \sum c_i a_i \big\|$. Thus,
since $\phi$ is contractive, its
restriction to $\spn \{ a_1 , \dots , a_m \}$ is a
$(1 - \delta \| \Gamma^{-1} \| )^{-1}$-isomorphism onto its image in
$B$. It follows that the
composition $\phi\circ\Gamma$ is a $\| \Gamma \|
(\| \Gamma^{-1} \|^{-1} - \delta )^{-1}$-isomorphism
onto its image in $B$. Since $B$ embeds into a matrix algebra of the same
rank, by Lemma 3.1 of \cite{EID} we reach the desired conclusion.
\end{proof}

For a function $f\in C(X)$ where $X$ is a compact Hausdorff space,
we denote by $d_f$ the pseudo-metric on $X$ given by
$$ d_f (x,y) = | f(x) - f(y) | $$
for all $x,y\in X$.

\begin{theorem}\label{T-comm}
Let $X$ be a compact metric space and $T : X\to X$ a homeomorphism.
Then for any $f\in C(X)_{\sa}$ the following are equivalent:
\begin{enumerate}
\item $ht(\alpha_T , \{ f \} ) > 0$,
\item there exists an entropy pair $(x,y)\in X\times X$ with $f(x) \neq f(y)$,
\item $h_{d_f} (T) > 0$,
\item there exist $K, d>0$, a sequence
$\{n_k\}_{k\in \mathbb{N}}$ in $\mathbb{N}$ going to $\infty$, and sets
$I_k\subseteq \{0, 1, \cdots, n_k-1\}$ of cardinality at least
$dn_k$ for each $k\in\mathbb{N}$ such that, for each $k\in\mathbb{N}$,
the linear map from $\ell^{I_k}_1$ to
$\spn \{ \alpha_T^i (f) : i\in I_k\}$
which sends the standard basis element of $\ell^{I_k}_1$
associated with $i\in I_k$ to $\alpha_T^i (f)$ is a $K$-isomorphism,
\item $f$ has an $\ell_1$ isomorphism set of positive density.
\end{enumerate}
\end{theorem}

\begin{proof}
(1)$\Rightarrow$(2)$\Rightarrow$(3). These implications follow
from the proofs of Theorem 4.3 and Lemma 4.2, respectively, in \cite{EID}.

(3)$\Rightarrow$(4)$\Rightarrow$(5). Here we can
apply the same arguments as in the proofs of the respective implications
(1)$\Rightarrow$(2)$\Rightarrow$(3) in Theorem~\ref{T-qs}.

(5)$\Rightarrow$(1). By assumption there exists a set $I\subseteq\mathbb{Z}$
of density greater than some $d>0$ and an isomorphism $\Gamma : \ell^I_1
\to\overline{\spn} \{ \alpha_T^i (f) : i\in I \}$ sending
the standard basis element of $\ell^I_1$ associated with
$i\in I$ to $\alpha_T^i (f)$. Let $0 < \delta < \| \Gamma^{-1} \|^{-1}$.
By Lemma~\ref{L-ell1} for every $n\in\mathbb{N}$ we have
\begin{gather*}
\log\rcp (  \{ \alpha_T^i (f) : i\in I \cap \{ -n , -n+1 , \dots , n \} \} ,
\delta ) \hspace*{30mm} \\
\hspace*{40mm} \geq | I \cap \{ -n , -n+1 , \dots , n \} | \cdot c
\cdot \| \Gamma \|^{-2} (\| \Gamma^{-1} \|^{-1} - \delta )^2
\end{gather*}
for some universal constant $c>0$. Since $I$ has density greater than $d$ and
$$ \rcp (\{ f , \alpha_T (f) , \dots , \alpha_T^{2n} (f) \} , \delta ) =
\rcp (\{ \alpha_T^{-n} (f)  , \alpha_T^{-n+1} (f) ,\cdots, \alpha_T^n(f)
\} , \delta ) , $$ we infer that
$$ \log\rcp (\{ f , \alpha_T (f) , \dots , \alpha_T^{2n} (f) \} ,
\delta )\geq d(2n+1)c \| \Gamma \|^{-2}
(\| \Gamma^{-1} \|^{-1} - \delta )^2 $$
for all sufficiently
large $n\in\mathbb{N}$. Hence $ht(\alpha_T , \{ f \}) > 0$.
\end{proof}

\begin{remark}\label{R-posentrcomm}
(i) By restricting to the $\alpha_T$-invariant unital
$C^*$-subalgebra of $C(X)$ generated by $f$ (which is separable and therefore
has metrizable pure state space) we need only
assume that $X$ is a compact Hausdorff space for the equivalences
(1)$\Leftrightarrow$(3)$\Leftrightarrow$(4)$\Leftrightarrow$(5)
in Theorem~\ref{T-comm}.

(ii) Since $h_{\topol}(\tilde{T}_{\alpha_T}) = \infty$ if and only
if $h_{\topol}(T) > 0$ (see \cite{GW} and \cite{EID}) if and only
if there is an entropy pair in $X\times X$ \cite{Disj},
Theorem~\ref{T-comm} shows that whenever $A$ is unital and
commutative and $Y\subseteq A_{\sa}$ generates $A$ as a unital
$C^*$-algebra (or, equivalently, separates the pure states of $A$)
we can require $a\in Y$ in conditions (2) and (3) in
Theorem~\ref{T-qs}.
\end{remark}

\begin{corollary}\label{C-Aff}
Let $A$ be a unital $C^*$-algebra and $\alpha$ an automorphism of
$A$. Let $\alpha'$ be the automorphism of $C(S(A))$ given by
$\alpha' (f)(\sigma ) = f(\sigma\circ\alpha )$ for all $f\in C(S(A))$ and
$\sigma\in S(A)$. Let $a\mapsto\bar{a}$ be the order isomorphism
from $A$ to the affine function space $\Aff (S(A)) \subseteq C(S(A))$ given by
$\bar{a} (\sigma ) = \sigma (a)$ for all $a\in A$ and $\sigma\in S(A)$.
Let $a\in A_{\sa}$. Then $ht(\alpha', \{ \bar{a} \} )>0$ implies
$ht(\alpha, \{ a \} )>0$.
\end{corollary}

\begin{proof}
The map $a\mapsto\bar{a}$ is an isomorphism of Banach spaces which
is isometric on $A_{\sa}$ and conjugates $\alpha$ to $\alpha'
\big| {}_{\Aff (S(A))}$, and so we obtain the conclusion from the
implication (1)$\Rightarrow$(4) in Theorem~\ref{T-comm} and an
appeal to Lemma~\ref{L-ell1}.
\end{proof}

\begin{theorem}\label{T-real}
The real topological Pinsker algebra is equal to the set of all
$f\in C(X)_{\sa}$ which do not have an $\ell_1$ isomorphism set of
positive density.
\end{theorem}

\begin{proof}
By Theorem 4.3 of \cite{EID} a function $f\in C(X)$ lies in the topological
Pinsker algebra if and only if $ht(\alpha_T ,\{ f \}) > 0$, and so an appeal 
to Theorem~\ref{T-comm} yields the result.
\end{proof}

Since $P_{X,T} = (P_{X,T})_{\sa} + i (P_{X,T})_{\sa}$ we obtain the
following corollary.

\begin{corollary}
The topological Pinsker algebra is equal to the
set of all $f\in C(X)$ whose real and imaginary
parts both lack an $\ell_1$ isomorphism set of positive density.
\end{corollary}

Since positive topological entropy implies the existence of an entropy pair
\cite{Disj}, the topological Pinsker algebra is equal to $C(X)$
precisely when $T$ has zero topological entropy and is equal to the 
scalars precisely when the system $(X,T)$ has completely positive
entropy (which means that every non-trivial factor
has positive topological entropy \cite{FPTE}).
We thus also have the following two corollaries.

\begin{corollary}\label{C-pe}
The homeomorphism $T$ has positive topological entropy if and only if
there is an $f\in C(X)$ with an $\ell_1$ isomorphism set of positive
density.
\end{corollary}

\begin{corollary}\label{C-cpe}
The system $(X,T)$ has completely positive entropy if and only if every
non-constant $f\in C(X)_{\sa}$ has an $\ell_1$ isomorphism set of
positive density.
\end{corollary}

\end{document}